\documentclass{amsart}
\usepackage{amsmath}
\usepackage{amssymb}
\usepackage{graphicx}
\input xy
\xyoption{all}
\newtheorem{Lemma}{Lemma}[section]
\newtheorem{Th}[Lemma]{Theorem}

\newenvironment{Proof}{{\sc Proof.}\ }{~\rule{1ex}{1ex}\vspace{0.2truecm}}
\newcommand{\Cal}[1]{{\mathcal #1}}

 \newcommand{\sign}{\operatorname{sign}}
    \newcommand{\Char}{\operatorname{Char}}

\begin{document}
    \title[Maximal subfields of a division algebra]{Maximal subfields of a division algebra}
    \author{Mai Hoang Bien}
\address{Mathematisch Instituut, Leiden Universiteit, Niels Bohrweg 1,2333 CA Leiden,The Netherlands.}
\curraddr{Dipartimento di Matematica, Universit\`{a} degli Studi di Padova, Via Trieste 63, 35121 Padova, Italy.}

\email{maihoangbien012@yahoo.com}
\thanks{The author would like to thank his supervisor Prof. H.W. Lenstra for the comments.}


   \maketitle
 \begin{abstract} Let $D$ be a division algebra over a field $F$. In this paper, we prove that there exist $a,b,x,y\in D^*=D\backslash\{0\}$ such that $F(ab-ba)$ and $F(xyx^{-1}y^{-1})$ are maximal subfields of $D$, which answers questions posted in \cite{Mah1}.  \end{abstract}

\section{Introduction}

Let $F$ be a field. A ring $D$ is called a {\it division algebra} over $F$ if the center $Z(D)=\{\,a\in D\mid ab=ba, \forall b\in D\,\}$ of $D$ is equal to $F$, $D$ is a finite dimensional vector space over $F$ and $D$ has neither proper left ideal nor proper right ideal. In other words, $D$ is a division ring with the center $F$ and $\dim_FD<\infty$. In some books and papers, $D$ is also called {\it centrally finite} \cite[Definition 14.1]{Lam}. A central simple algebra over $F$ is an algebra isomorphic to $M_n(D)$ for some positive integer $n$ and division algebra $D$ over $F$. For any central simple algebra $A$ over $F$, $\sqrt{\dim_FA}$ is said to be {\it degree} of $A$.

For any division algebra $D$ over $F$, it is well known from Kothe's Theorem that there exists a maximal subfield $K$ of $D$ such that the extension of fields $K/F$ is separable \cite[Th. 15.12]{Lam}. In \cite[Theorem 7]{AAM}, authors proved that for any separable extension  of fields  $K/F$ in $D$, there exists an element $c\in [D,D]$, the group of additive commutators of $(D,+)$, such that $K=F(c)$ unless $\Char (F)=[K:F]=2$ and $4$ does not divide the degree of $D$. Hence, if $K$ is a maximal subfield of $D$ which is separable over $F$ then there exists $c\in [D, D]$ such that $K=F(c)$. In particular, there exists a maximal subfield of $D$ such that it is of the form $F(c)$ for some element $c$ in $[D,D]$. We have a natural question: is it true that there exists a commutator $ab-ba\in [D,D]$ such that $F(ab-ba)$ is a maximal subfield of $D$ (see \cite[Problem 28]{Mah1})?  Almost similarly, if $K/F$ is a separable extension of fields in $D$ then there exists an element $d\in D'=[D^*,D^*]$, the group of multiplicative commutators of $D^*=D\backslash \{0\}$, such that $K=F(d)$ (see \cite[Theorem 2.26]{Mah1}). Again, the author asked whether $F(xyx^{-1}y^{-1})$ is a maximal subfield of $D$ for some $x,y\in D^*$ (see \cite[Problem 29]{Mah1}). 

The goal of this paper is to answer in the affirmative for both questions. The main tools used in this paper are generalized rational identities over a central simple algebra. Readers can find their definitions and notaions in detail in \cite{BMM} and \cite{ Rowen}. 

\section{Results} 

Let $R$ be a ring. Recall that an element $a$ of $R$ is called {\it algebraic of degree $n$} over a subring $S$ of $R$  if there exists a polynomial $f(x)$  of degree $n$  over $S$ such that $f(a)=0$ and there is no polynomial of degree less than $n$ vanishing on $a$. In general, $f(x)$ is not necessary unique and irreducible even if $S$ is a field. For example, the matrix $A=\left( {\begin{array}{*{20}{c}}
1&0\\
0&2
\end{array}} \right)\in M_2(F), $ where $F$ is a field, satisfies the polynomial $f(x)=(x-1)(x-2)$. Since $A\notin F$, $2$ is the smallest degree of all the polynomials vanishing on $A$. 

Recall that a {\it generalized rational expression} over $R$ is an expression contructed from $R$ and a set of noncommutative inderteminates  using addition, substraction, multiplication and division. A generalized rational expression $f$ over $R$ is called a {\it generalized rational identity} if it vanishes on all permissible substitutions from $R$. In this case, one says {\it $R$ satisfies $f$}. We consider the following example which is important in this paper. Given a positive integer $n$ and $n+1$ noncommutative indeterminates $x,y_1,\cdots, y_n$, put $$g_n(x,y_1,y_2,\cdots, y_n)=\sum\limits_{\delta  \in {S_{n + 1}}} {\sign(\delta ){x^{\delta (0)}}{y_1}{x^{\delta (1)}}{y_2}{x^{\delta (2)}} \ldots {y_n}{x^{\delta (n)}}}, $$ where $S_{n+1}$ is the symmetric group of $\{\,0,1,\cdots, n\,\}$ and $\sign(\delta)$ is the sign of permutation $\delta$. This is a generalized rational expression defined in \cite{BMM} to connect an algebraic element of degree $n$ and a polynomial of $n+1$ indeterminates.  
\begin{Lemma}\label{2.2}Let $F$ be a field and $A$ be a central simple algebra over $F$. For any element $a\in A$, the following conditions are equivalent.
\begin{enumerate}
\item The element $a$ is algebraic over $F$ of degree less than $n$.
\item $g_n(a,r_1,r_2,\cdots, r_n)=0$ for any $r_1, r_2,\cdots, r_n\in A$.
\end{enumerate} 
\end{Lemma}
\begin{Proof} This is a corollary of \cite[Corollary 2.3.8]{BMM}.
\end{Proof}

In particular, a central simple algebra of degree $m$ satisfies the expression $g_m$ since every central simple algebra of degree $m$ over a field $F$ can be considered as a $F$-subalgebra of the ring $M_m(F)$ and elements of $M_m(F)$ are algebraic of degree less than $m$ over $F$. In other words, $g_m$ is a generalized rational indentity of any central simple algebra of degree $m$. 

For any central simple algebra $A$, denote $\Cal G(A)$ the set of all generalized rational identities of $A$. Then $\Cal G(A)\ne \emptyset$ because $g_m\in \Cal G(A)$. The following theorem gives us a relation between the set of all generalized rational identities of a central simple algebra and the ring of matrices over a field.

\begin{Th}\label{2.3}\cite[Theorem 11]{Ami1} Let $F$ be an infinite field and $A$ be a central simple algebra of degree $n$ over $F$. Assume that $L$ is an extension field of $F$. Then $\Cal G(A)=\Cal G(M_n(F))=\Cal G(M_n(L))$.
\end{Th}

Now we are going to prove the main results of this paper. The following lemma is basic.

\begin{Lemma}\label{3.0} Let  $D$ be a division algebra of degree $n$ over a field $F$. Assume that $K$ is a subfield of $D$ containing $F$. Then $\dim_FK\le n$. The quality holds if and only if $K$ is a maximal sufield of $D$.
\end{Lemma}
\begin{Proof} See \cite[Corollary 15.6 and Proposition 15.7]{Lam}
\end{Proof}

\begin{Lemma} \label{3.1}Let $F$ be an infinite field and $n\ge 2$ be an integer. There exist two matrices $A,B\in M_n(F)$ such that the commutator $ABA^{-1}B^{-1}$ is an algebraic element of degree $n$ over $F$. 
\end{Lemma}
\begin{Proof} Put 

$A=\left( {\begin{array}{*{20}{c}}
0&0& \cdots &0&{{a_1}}\\
1&0& \cdots &0&{{a_2}}\\
 \vdots & \vdots & \ddots & \vdots & \vdots \\
0&0&1&0&{{a_{n - 1}}}\\
0&0&0&1&0
\end{array}} \right)$ and $B=\left( {\begin{array}{*{20}{c}}
{{b_1}}&0& \cdots &0&0\\
0&{{b_2}}& \cdots &0&0\\
 \vdots & \vdots & \ddots & \vdots & \vdots \\
0&0&0&{{b_{n - 1}}}&0\\
0&0&0&0&{{b_n}}
\end{array}} \right),$ where $a_i,b_j\ne 0$. One has $ABA^{-1}B^{-1}=\left( {\begin{array}{*{20}{c}}
{{b_n}b_1^{ - 1}}&0& \cdots &0&0\\
*&{{b_1}b_2^{ - 1}}& \cdots &0&0\\
 \vdots & \vdots & \ddots & \vdots & \vdots \\
*&*&*&{{b_{n - 2}}b_{n - 1}^{ - 1}}&0\\
*&*&*&*&{{b_{n - 1}}b_n^{ - 1}}
\end{array}} \right)$. If we choose $b_nb_1^{-1}, b_1b^{-1}_2,\cdots, b_{n-1}b^{-1}$ all distinct (it is possible since $F$ is infinite), then the characteristic polynomial of $ABA^{-1}B^{-1}$ is a polynomial of smallest degree which vanishes on $ABA^{-1}B^{-1}$. That is, $ABA^{-1}B^{-1}$ is an algebraic element of degree $n$ over $F$.
\end{Proof}

The following theorem answers Problem 29 in \cite[Page 83]{Mah1}.
\begin{Th}\label{3.2} Let $D$ be a central division algebra over a field $F$. There exist $x,y\in D^*$ such that $F(xyx^{-1}y^{-1})$ is a maximal subfield of $D$.
\end{Th}

\begin{Proof} If $F$ is finite then $D$ is also finite, so that there is nothing to prove. Suppose that $F$ is infinite and $D$ is of degree $n$ over $F$. By Lemma~\ref{3.0}, it suffices to show that there exist $x,y\in D^*$ such that $\dim_FF(xyx^{-1}y^{-1})\ge n$. Indeed, put $\ell=\max\{\,\dim_FF(xyx^{-1}y^{-1})\mid x,y\in D^*\,\}.$ Then from Lemma~\ref{3.0}, $$g_\ell(rsr^{-1}s^{-1}, r_1, r_2,\cdots ,r_\ell )=0$$ for any $r_1,r_2,\cdots , r_\ell\in D$ and $r,s\in D^*$. Hence, $g_\ell(xyx^{-1}y^{-1}, y_1,y_2,\cdots, y_\ell)$ is a generalized rational idenity of $D$, so that, by Lemma~\ref{2.3}, $g_\ell(xy-yx, y_1,y_2,\cdots, y_\ell)$ is also a generalized rational idenity of $M_n(F)$. Since $g_\ell(ABA^{-1}B^{-1}, r_1, r_2,\cdots ,r_\ell )=0, $ for any $ r_i\in M_n(F)$ and $A,B$ are chosen in Lemma~\ref{3.1}. Therefore $n\le \ell$ because  Lemma~\ref{2.2} and $AB-BA$ is an algebraic element of degree $n$.\end{Proof}

\begin{Lemma} \label{3.3}Let $F$ be an infinite field and $n>2$ be an integer. There exist two matrices $A,B\in M_n(F)$ such that $AB-BA$ is an algebraic element of degree $n$ over $F$. 
\end{Lemma}

\begin{Proof} Put 

$A=\left( {\begin{array}{*{20}{c}}
0&0& \cdots &0&{{a_1}}\\
1&0& \cdots &0&{{a_2}}\\
 \vdots & \vdots & \ddots & \vdots & \vdots \\
0&0&1&0&{{a_{n - 1}}}\\
0&0&0&1&0
\end{array}} \right)$ and $B=\left( {\begin{array}{*{20}{c}}
0&{{b_1}}&0& \cdots &0&0\\
0&0&{{b_2}}& \cdots &0&0\\
0&0&0& \cdots &0&0\\
 \vdots & \vdots & \vdots & \ddots & \vdots & \vdots \\
0&0&0& \cdots &0&{{b_{n - 1}}}\\
0&0&0& \cdots &0&0
\end{array}} \right)$. One has $AB-BA=\left( {\begin{array}{*{20}{c}}
{{b_1}}&*& \cdots &*&*\\
0&{{b_1} - {b_2}}& \cdots &*&*\\
 \vdots & \vdots & \ddots & \vdots & \vdots \\
0&0& \cdots &{{b_{n - 2}} - {b_{n - 1}}}&*\\
0&0& \cdots &0&{{b_{n - 1}}}
\end{array}} \right)$. Since $F$ is infinite, we can choose $b_1,b_2,\cdots, b_{n-1}\in F$ such that $b_1,b_1-b_2,\cdots, b_{n-2}-b_{n-1}, b_{n-1}$ all distinct. Hence, the characteristic polynomial of $AB-BA$ is a polynomial of smallest degree vanishing on $AB-BA$. Therefore, $AB-BA$ is an algebraic element of degree $n$ over $F$.
\end{Proof}

Almost similar to the proof of Theorem~\ref{3.2}, we have the following theorem, which answers Problem 28 in \cite[Page 83]{Mah1}.
\begin{Th}\label{3.4} Let $D$ be a central division algebra over a field $F$. There exist $x,y\in D$ such that $F(xy-yx)$ is a maximal subfield of $D$.
\end{Th}

\begin{Proof} If $F$ is finite then $D$ is also finite, so that there is nothing to prove. Suppose that $F$ is infinite and $D$ is of degree $n$. By Lemma~\ref{3.0}, it suffices to show that there exist $x,y\in D$ such that $\dim_FF(xy-yx)\ge n$. Indeed, if $n=2$, by \cite[Corollary 13.5]{Lam}, then there exist $x,y\in D$ such that $xy-yx\notin F$, which implies $F(xy-yx)=2=n$. Assume that $n>2$. Then put $\ell=\max\{\,\dim_FF(xy-yx)\mid x,y\in D\,\}.$ By Lemma~\ref{2.2}, $$g_\ell(rs-sr, r_1, r_2,\cdots ,r_\ell )=0$$ for any $r_1,r_2,\cdots , r_\ell\in D$ and $r,s\in D^*$. It follows $g_\ell(xy-yx, y_1,y_2,\cdots, y_\ell)$ is a generalized rational idenity of $D$. From Lemma~\ref{2.3}, $g_\ell(xy-yx, y_1,y_2,\cdots, y_\ell)$ is also a generalized rational idenity of $M_n(F)$. But because there exist $A, B\in M_n(F)$ such that $AB-BA$ is algebraic of degree $n$ (Lemma~\ref{3.1}), one has $$g_\ell(AB-BA, r_1, r_2,\cdots ,r_\ell )=0$$ for any $r_i\in M_n(F)$. Therefore, by Lemma~\ref{2.2}, $n\le \ell$.\end{Proof}


\begin{thebibliography}{}
\bibitem{AAM} S. Akbari,  M. Arian-Nejad, M. L. Mehrabadi, On additive commutator groups in division rings, {\it Results Math.}, {\bf 33} (1-2), 9–21, 1998.
\bibitem{Ami1} S. A. Amitsur, Rational identities and applications to algebra and geometry, {\it J. Algebra} {\bf 3}, 304--359, 1966.

\bibitem{BMM} K. I. Beidar, W. S. Martindale 3rd and A. V. Mikhalev, {\it Rings with Generalized Identities}, Marcel Dekker, Inc., New York- Basel-Hong Kong, 1996.
\bibitem{Lam} T. Y. Lam, {\it A first course in noncommutative rings}, MGT 131, Springer, 1991.
\bibitem{Mah1} M. Mahdavi-Hezavehi, Commutators in division rings revisited. {\it Bull. Iranian Math. Soc}, {\bf 26}(3): 7--88, 2000.
\bibitem{Rowen} L. H. Rowen, {\it Polynomial identities in ring theory}, Academic Press, Inc., New York, 1980.
 \end{thebibliography}
\end{document}